\documentclass[a4paper,12pt]{article}

\usepackage{amsfonts}
\usepackage{amssymb}
\usepackage{latexsym}
\usepackage{amsthm}
\usepackage[dvips]{graphicx}

\usepackage{fancyhdr}

\begin{document}

\theoremstyle{plain}       \newtheorem*{mainthm}{Theorem 1}
\theoremstyle{plain}       \newtheorem*{cor1}{Corollary 1}
\theoremstyle{plain}       \newtheorem*{cor2}{Corollary 2}
\theoremstyle{plain}       \newtheorem*{thm2}{Theorem 2}
\theoremstyle{plain}       \newtheorem{thm}{Theorem}[section]
\theoremstyle{definition}  \newtheorem{defn}[thm]{Definition}
\theoremstyle{plain}       \newtheorem{lemma}[thm]{Lemma}
\theoremstyle{plain}       \newtheorem{prop}[thm]{Proposition}

\title{Consequences of Ergodic Frame Flow for Rank Rigidity in Negative Curvature}
\author{David Constantine\footnote{Supported by NSF Graduate Research Fellowship}   
\footnote{\textsc{Department of Mathematics, University of Michigan, Ann Arbor, MI 48103 
U.S.A.}                  \textit{email:} \texttt{constand@umich.edu}}}
\date{December 19, 2006}
\maketitle


\begin{abstract}
This paper presents a rank rigidity result for negatively curved spaces.  Let $M$ be a 
compact manifold with negative sectional curvature and suppose that along every geodesic 
in $M$ there is a parallel vector field making curvature $-a^2$ with the geodesic 
direction.  We prove that $M$ has constant curvature equal to $-a^2$ if $M$ is odd 
dimensional, or if $M$ is even dimensional and has sectional curvature pinched as 
follows: $-\Lambda^2 < K < -\lambda^2$ where $\lambda/\Lambda > .93$.  When $a$ is 
extremal, i.e. $-a^2$ is the curvature minimum or maximum for the manifold, this result 
is analogous to rank rigidity results in various other curvature settings where higher 
rank implies that the space is locally symmetric.  In particular, this result is the 
first positive result for lower rank (i.e. when $-a^2$ is minimal), and in the upper rank 
case gives a shorter proof of the hyperbolic rank rigidity theorem of Hamenst\"{a}dt, 
subject to the pinching condition in even dimension.  We also present a rigidity result 
using only an assumption on maximal Lyapunov exponents in direct analogy with work done 
by Connell.  Our proof of the main theorem uses the ergodic theory of the frame flow 
developed by Brin and others - in particular the transitivity group associated to this 
flow.
\end{abstract}


\section{Introduction}
Rank rigidity was first proved in the higher Euclidean rank setting by Ballmann 
\cite{Ballman} and, using different methods, by Burns and Spatzier \cite{Burns-Spatzier}. 
A manifold is said to have higher Euclidean rank if a parallel normal Jacobi field can be 
found along every geodesic.  Ballmann and Burns-Spatzier proved that if an irreducible, 
compact, nonpositively curved manifold has higher Euclidean rank, then it is locally 
symmetric.  Ballmann's proof works for finite volume as well and the most general version 
of this theorem is due to Eberlein and Heber, who prove it under only a dynamical 
condition on the isometry group of $M$'s universal cover \cite{Eberlein-Heber}.  
Hamenst\"{a}dt showed that a compact manifold with curvature bounded above by -1 is 
locally symmetric if along every geodesic there is a Jacobi field making curvature -1 
with the geodesic direction \cite{Hamenstadt}.  She called this situation higher 
hyperbolic rank.  Shankar, Spatzier and Wilking extended rank rigidity into positive 
curvature by defining spherical rank.  A manifold with curvature bounded above by 1 is 
said to have higher spherical rank if every geodesic has a conjugate point at $\pi$, or 
equivalently, a parallel vector field making curvature 1 with the geodesic direction.  
They proved that a complete manifold with higher spherical rank is a compact, rank one 
locally symmetric space \cite{sphericalrank}.

These results settle many rank rigidity questions, but leave questions about other 
curvature settings open (see \cite{sphericalrank} for an excellent overview).  In this 
paper we prove the following theorem, which can be applied to various settings in 
negative curvature.

\begin{mainthm} Let $M$ be a compact, negatively curved manifold.  Suppose that along 
every geodesic in $M$ there exists a parallel vector field making sectional curvature 
$-a^2$ with the geodesic direction.  If $M$ is odd dimensional, or if $M$ is even 
dimensional and satisfies the sectional curvature pinching condition $-\Lambda^2 < K < 
-\lambda^2$ with $\lambda/\Lambda > .93$ then $M$ has constant negative curvature equal 
to $-a^2$.
\end{mainthm}

Note that, unlike previous rank rigidity results, Theorem 1 allows for situations where 
the distinguished curvature $-a^2$ is not extremal.  However, the cases where $-a^2$ is 
extremal are of particular importance and in these situations the extremality of the 
distinguished curvature $-a^2$ allows the hypotheses of our theorem to be weakened, as 
demonstrated in section \ref{sec:parallel} of this paper.  The folowing two results are 
then easy corollaries of Theorem 1:

\begin{cor1}
Let $M$ be a compact manifold with sectional curvature $-1 \leq K < 0$.  Suppose that 
along every geodesic in $M$ there exists a Jacobi field making sectional curvature $-1$ 
with the geodesic direction.  If $M$ is odd dimensional, or if $M$ is even dimensional 
and satisfies the sectional curvature pinching condition $-1 \leq K < -\lambda^2$ with 
$\lambda > .93$ then $M$ is hyperbolic.
\end{cor1}

\begin{cor2} \emph{(compare with Hamenst\"{a}dt \cite{Hamenstadt})}
Let $M$ be a compact manifold with sectional curvature bounded above by $-1$.  Suppose 
that along every geodesic in $M$ there exists a Jacobi field making sectional curvature 
$-1$ with the geodesic direction.  If $M$ is odd dimensional, or if $M$ is even 
dimensional and satisfies the sectional curvature pinching condition $-(1/.93)^2 < K \leq 
-1$ then $M$ is hyperbolic.
\end{cor2}

In Corollary 1, $-1$ is the curvature minimum for $M$ and we obtain a new rank rigidity 
result analogous to those described above.  This is the first positive result for lower 
rank, i.e. when the distinguished curvature value is the lower curvature bound (see 
section \ref{sec:conclusion} for more discussion).  In Corollary 2, $-1$ is the curvature 
maximum for $M$ and we obtain a shorter proof of Hamenst\"{a}dt's result, under an added 
pinching constraint in even dimension.

In \cite{Connell}, Connell showed that rank rigidity results can be obtained using only a 
dynamical assumption on the geodesic flow, namely an assumption on the Lyapunov exponents 
at a full measure set of unit tangent vectors.  His paper deals with the upper rank 
situations treated by Ballmann, Burns-Spatzier and Hamenst\"{a}dt.  He proves that having 
the minimal Lyapunov exponent allowed by the curvature restrictions attained at a full 
measure set of unit tangent vectors is sufficient to apply the results of Ballman and 
Burns-Spatzier or Hamenst\"{a}dt.  In the lower rank setting of this paper, this 
viewpoint translates into

\begin{thm2}
Let $M$ be a compact manifold with sectional curvature $-a^2 \leq K < 0$, where $a > 0$.  
Suppose that for a full (Liouville) measure set of unit tangent vectors $v$ on $M$ the 
maximal Lyapunov exponent at $v$ is $a$, the maximum allowed by the curvature 
restriction.  If $M$ is odd dimensional, or if $M$ is even dimensional and satisfies the 
sectional curvature pinching condition $-a^2 \leq K < -\lambda^2$ with $\lambda/a > .93$ 
then $M$ is of constant curvature $-a^2$.
\end{thm2}

\noindent The adaptation of Connell's arguments for this setting is discussed in section 
\ref{sec:Lyap}.

The proof of Theorem 1 relies on dynamical properties of the geodesic and frame flows on 
negatively curved manifolds.  We rely heavily on Brin's work on frame flows (see 
\cite{Brin-survey} for a survey), and results on the ergodicity of these flows due to 
Brin and Gromov in odd dimension \cite{Brin-Gromov} and to Brin and Karcher in even 
dimension \cite{Brin-Karcher}.  These results are summarized in section 
\ref{sec:background}.  In particular, we utilize the transitivity group $H_v$, defined 
for any vector $v$ in the unit tangent bundle of $M$, which acts on $v^{\perp} \subset 
T^1M$.  Essentially, elements of $H_v$ correspond to parallel translations around ideal 
polygons in $M$'s universal cover.  Brin shows that this group is the structure group for 
the ergodic components of the frame flow (see e.g \cite{Brin-survey} or 
\cite{Brin-gpext}).  In section \ref{sec:trans} we show that, subject to suitable 
recursion properties on these ideal polygons, $H_v$ preserves the parallel fields that 
make curvature $-a^2$ with the geodesic defined by $v$.  We then show that these 
recursion properties are generic and that the elements of $H_v$ vary continuously with 
the choice of ideal polygon.  Thus, all of $H_v$ respects the distinguished fields.  
Finally, we apply results of Brin-Gromov and Brin-Karcher on the ergodicity of the 
2-frame flow which imply that $H_v$ acts transitively on $v^\perp$ and conclude that the 
curvature of $M$ is constant.

I would like to thank Chris Connell for discussions helpful with the arguments in section 
\ref{sec:parallel} of this paper, Jeffrey Rauch for the proof of Lemma \ref{Rauchlemma}, 
and Ben Schmidt for helpful comments on this paper.  In particular, special thanks are 
due to my advisor, Ralf Spatzier, for suggesting this problem, for help with several 
pieces of the argument and for helpful comments on this paper.


\section{Notation and background} \label{sec:background}

Let us begin by fixing some notation and stating the results we will need.  Let $M$ be a 
compact Riemannian manifold with negative sectional curvature and let $\tilde{M}$ be its 
universal cover.  Denote by $T^1M$ and $T^1\tilde{M}$ the unit tangent bundles to $M$ and 
$\tilde{M}$, respectively.  We will denote by $g_t$ the geodesic flow on either of these 
spaces, and by $F_t$ the frame flow on the Stiefel manifold $St_kM$, the space of ordered 
orthonormal $k$-frames on $M$.  $St_kM$ is a fiber bundle over $T^1M$ with the group 
$SO(n-1)$ acting on the right; $St_nM$ is a principal bundle with $SO(n-1)$ as structure 
group.  There are standard measures on these spaces, namely Liouville measure on $T^1M$ 
and $T^1\tilde{M}$, and on $St_kM$ the product measure of Liouville measure and the 
measure on the fibers inherited from the Haar measure on $SO(n-1)$.  Unless otherwise 
specified, these will be the measures used in all that follows.  Let $\gamma_v (t)$ 
denote the geodesic in $M$ or $\tilde{M}$ with velocity $v$ at time 0.  We will denote by 
$w_v (t)$ a parallel normal vector field along $\gamma_v (t)$ making the distinguished 
curvature $-a^2$ with $\dot{\gamma}_v (t)$.  Finally, $\langle \cdot, \cdot \rangle$ will 
denote the Riemannian inner product, $R(\cdot, \cdot) \cdot$  will denote the curvature 
tensor and $K(\cdot, \cdot)$ will denote sectional curvature.

The geodesic flow gives rise to stable and unstable foliations $W^s_g$ and $W^u_g$ of 
$T^1M$ or $T^1\tilde{M}$.  These foliations are absolutely continuous, and the geodesic 
flow for such an $M$ is ergodic (proved by Anosov, see Brin's appendix to 
\cite{Brin-ergodicity}).  The frame flow also gives rise to stable and unstable 
foliations $W^s_F$ and $W^u_F$ as shown by Brin \cite{Brin-toptrans}.  These foliations 
allow Brin to define the \emph{transitivity group} in the following way.  If $v$ and $v'$ 
are on the same leaf of $W^s_g$ then for every $n$-frame $\alpha$ above $v$ there is a 
unique $n$-frame $\alpha'$ above $v'$ such that $\alpha$ and $\alpha'$ belong to the same 
leaf of $W^s_F$.  In particular, the distance between $F_t(\alpha)$ and $F_t(\alpha')$ 
approaches 0 as $t \rightarrow \infty$ (this is how one determines that $\alpha'$ belongs 
to the leaf $W^s_F(\alpha)$).  Let $p(v, v')$ be the map from the fiber of $St_nM$ over 
$v$ to the fiber over $v'$ that takes each $\alpha$ to the corresponding $\alpha'$ over 
$v'$.  Note $p(v, v')$ corresponds to a unique isometry between $v^{\perp}$ and 
$v'^{\perp}$.  Once defined on $n$-frames, $p(v, v')$ acts on all $k$-frames; the action 
on 2-frames will be what we use in this paper and we will abuse notation by using $p(v, 
v')$ to denote this restricted action.  One can think of $p(v, v')(\alpha)$ as the result 
of parallel transporting $\alpha$ along $\gamma_v (t)$ out `to the boundary at infinity 
of $\tilde{M}$' and then back to $v'$ along $\gamma_{v'} (t)$.  If $v'$ and $v$ belong to 
the same leaf of $W^u_g$ there is similarly an isometry corresponding to parallel 
translation to the boundary at infinity along $\gamma_{-v}$ and back along 
$\gamma_{-v'}$.  This defines the unstable leaves for the frame flow foliation and we 
will also denote this isometry by $p(v, v')$.  Brin (see \cite{Brin-survey} Defn. 4.4) 
then defines the transitivity group at $v$ as follows:

\begin{defn}
Given any sequence $s=\{v_0, v_1, \ldots , v_k\}$ with $v_0=v_k=v$ such that each pair 
$\{v_i, v_{i+1}\}$ lies on the same leaf of $W^s_g$ or $W^u_g$ we have an isomorphism of 
$v^{\perp}$ given by \[I(s) = \prod_{i=0}^{k-1} p(v_i, v_{i+1}).\]
The closure of the set of all such isometries is denoted by $H_v$ and is called the 
transitivity group.
\end{defn}

\begin{figure}
\centering
\includegraphics[height=.45\textheight, width=.6\textwidth]{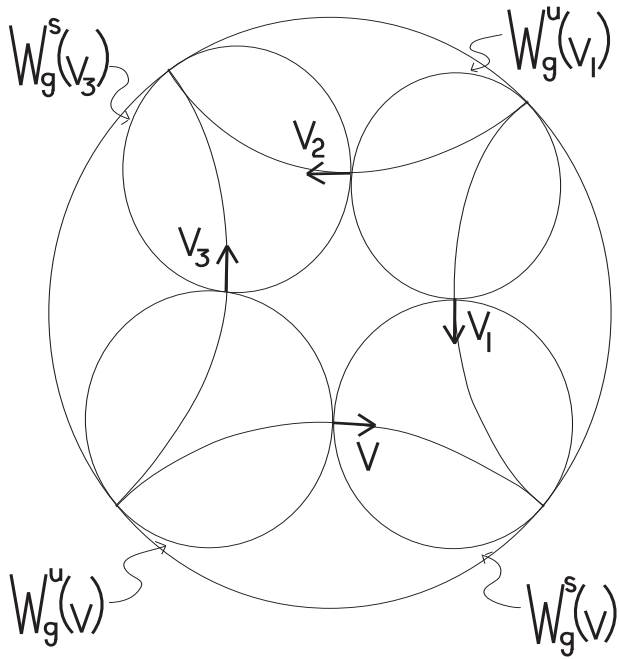}
\caption{`equilateral' ideal rectangle} \label{fig:equi}
\end{figure}

The idea of the transitivity group is that it is generated by isometries coming from 
parallel translation around ideal polygons in $\tilde{M}$ with an even number of sides, 
such as the one shown in figure \ref{fig:equi}.  Note that here only `equilateral' 
polygons are allowed; in figure \ref{fig:equi} only rectangles for which $W^u_g(v_1)$ is 
tangent to $W^s_g(v_3)$ are permitted.  We will later find it useful to allow general 
ideal polygons.

The definition of this group arises in Brin's analysis of the ergodic components of the 
frame flow.  He shows in \cite{Brin-gpext} that the ergodic components are subbundles of 
$St_kM$ with structure group a closed subgroup of $SO(n-1)$ (see also \cite{Brin-survey} 
section 5 for an overview).  In addition, his proof demonstrates that the structure group 
for the ergodic component is the transitivity group (see \cite{Brin-survey} Remark 2 or 
\cite{Brin-gpext} Proposition 2).  This explicit geometric description of the ergodic 
components is the central tool used in our proof.

We use two results on the ergodicity of the 2-frame flow in our proof.

\begin{thm} \label{thm:BG}
\emph{(Brin-Gromov \cite{Brin-Gromov} Proposition 4.3)} If $M$ has negative sectional 
curvature and odd dimension then the 2-frame flow is ergodic.
\end{thm}

\begin{thm} \label{thm:BK}
\emph{(Brin-Karcher \cite{Brin-Karcher})} If $M$ has sectional curvature satisfying 
$-\Lambda^2 < K < -\lambda^2$ with $\lambda / \Lambda > .93$ then the 2-frame flow is 
ergodic.
\end{thm}

Theorem \ref{thm:BK} is not directly stated as above in \cite{Brin-Karcher}, rather it 
follows from remarks made in section 2 of that paper together with Proposition 2.9 and 
the extensive estimates carried out in the later sections.  Note that since the 2-frame 
flow preserves the parallel fields making curvature $-a^2$, ergodicity of this flow alone 
seems to indicate that the manifold has constant curvature.  However, since the subset of 
$St_2M$ given by these distinguished fields may, a priori, have zero measure, the result 
does not follow directly from ergodicity.  Instead, we must use the precise description 
of the ergodic components given by the transitivity group.


\section{The transitivity group and distinguished vector fields} \label{sec:trans}

As noted in the Introduction, the description of the ergodic components in terms of the 
transitivity group is crucial.  In this section we investigate how the distinguished 
vector fields $w_v(t)$ along $\gamma_v(t)$ behave under the action of the transitivity 
group and use the results to prove Theorem 1.

To obtain results we will need to assume some dynamical properties of the ideal polygon 
that produces a given element of $H_v$.  For example, consider the ideal rectangle 
defined by $v$, $v_1$, $v_2$ and $v_3$ as pictured in figure \ref{fig:equi}.  Note that 
the distinguished vector field $w_v(t)$ which makes constant curvature $-a^2$ with 
$\gamma_v(t)$ corresponds uniquely to a parallel normal vector field $P(t)$ along 
$\gamma_{v_1}(t)$, such that $v_1$ and $P(0)$ make a 2-frame in the leaf of the stable 
foliation containing the 2-frame $\{v, w_v(0)\}$, that is, $\{v_1, P(0)\} = p(v, v_1)\{v, 
w_v(0)\}$.  By continuity of the sectional curvature, $K(P(t), \dot{\gamma}_{v_1}(t)) 
\rightarrow -a^2$ as $t \rightarrow \infty$.

\begin{lemma}
Suppose $\gamma(t)$ is a recurrent geodesic with a parallel normal field $P(t)$ along it 
such that $K(P(t), \dot{\gamma}(t)) \rightarrow -a^2$ as $t \rightarrow \infty$.  Then 
$K(P(t), \dot{\gamma}(t)) \equiv -a^2$ for all $t$.
\end{lemma}

\begin{proof}
Since $\gamma(t)$ is recurrent we can take an increasing sequence $\{ t_k\} $ tending to 
infinity such that $\dot{\gamma}(t_k)$ approaches $\dot{\gamma}(0)$.  Since the parallel 
field $P(t)$ has constant norm and the set of vectors in $\dot{\gamma}_v(t)^{\perp}$ with 
this norm is compact, we can, by passing to a subsequence, assume that $P(t_k)$ has a 
limit $G(0)$.  Extend $G(0)$ to a parallel vector field $G(t)$ along $\gamma(t)$.

By construction, $K(G(0), \dot{\gamma}(0)) = \lim_{k \rightarrow \infty} K(P(t_k), 
\dot{\gamma}(t_k)) = -a^2$.  In addition, for any real number $T$, the recurrence 
$\dot{\gamma}(t_k) \rightarrow \dot{\gamma}(0)$ implies recurrence $\dot{\gamma}(t_k+T) 
\rightarrow \dot{\gamma}(T)$.  By continuity of the frame flow, we get that $P(t_k+T) 
\rightarrow G(T)$ for the vector field $G$ defined above.  Thus $G(t)$ makes curvature 
$-a^2$ with $\dot{\gamma}(t)$ for any time $t$.

We can repeat the same argument as above, letting $G(t)$ recur along the same sequence of 
times to produce $G_1(t)$, and likewise $G_i(t)$ recur to produce $G_{i+1}(t)$, forming a 
sequence of fields all making curvature identically $-a^2$ with the geodesic direction.  
Now, observe that $G(0) = P(0)\cdot g$ for some $g\in SO(n-1)$.  Note here that $g$ is 
not well defined by looking at $P$ and $G$ alone, but will be well defined if we consider 
$n$-frame orbits with second vector $P$ recurring to $n$-frames with second vector 
$G(0)$; this is the $g$ we utilize.  By construction and the fact that the $SO(n-1)$ 
action commutes with parallel translation, $G_i(0)=P(0)\cdot g^{i+1}$.  $SO(n-1)$ is 
compact, so the $\{g^i\}$ have convergent subsequences.  In addition, since the terms of
this sequence are all iterates of a single element, we can, by adjusting terms of such a 
subsequence by suitable negative powers of $g$, have the subsequence converge to the 
identity.  Choose a subsequence $\{i_j\}$ such that $g^{i_j+1} \rightarrow id$ as $j 
\rightarrow \infty$.  These $G_{i_j}(t)$ approach our original field $P(t)$ showing that 
$P$ makes constant curvature $-a^2$ with $\dot{\gamma}$ as well.
\end{proof}

Consider the situation depicted in figure 1.  Lemma 3.1 shows that, when $\gamma_{v_1}$ 
is recurrent in forward time, the map $p(v, v_1)$ preserves the distinguished vector 
fields in the sense that it sends a vector from one such field, $w_v(0)$, to a vector 
from another such field along $\gamma_{v_1}$.  Thus, if in figure 1 we have that 
$\gamma_{v_1}$ and $\gamma_{v_3}$ are recurrent in positive time and $\gamma_v$ and 
$\gamma_{v_2}$ are recurrent in negative time, then the element of $H_v$ given by 
parallel translation around this ideal rectangle will map $w_v(0)$ to another element of 
$v^{\perp}$ which is in a parallel field along $\gamma_v$ making curvature $-a^2$.  If 
these sort of recurrence properties held for all `equilateral' ideal polygons based at 
$v$ we would have that the transitivity group preserves the distinguished vector fields.  
We cannot assure that these recurrence properties are always present, but ergodicity of 
the geodesic flow on $M$ indicates that they will be present almost all the time.  We now 
work out the details of this.

First, the ergodicity of the geodesic flow implies that there is a full measure set of 
vectors $v$ in $T^1M$ which have dense forward and backward orbits under the geodesic 
flow.  Choosing $v$ from this set implies that $\gamma_{v_1}$ will be recurrent in 
positive time and $\gamma_v$ will be recurrent in negative time.  For ideal rectangles, 
this leaves only the positive time recurrence of $\gamma_{v_3}$ and the negative time 
recurrence of $\gamma_{v_2}$ lacking.  It is convenient at this time to extend the 
definition of the transitivity group.

Consider the situation depicted in figure 2.  Here the unit tangent vectors $v$, $v_1$, 
$v_2$ and $v_3$ describe an ideal rectangle in $T^1\tilde M$.  Each pair $\{v, v_1\}$, 
$\{v_1, v_2\}$, $\{v_3, v\}$ lies on a leaf of $W^s_g$ or $W^u_g$ and we let $T\in 
\mathbb{R}$ be the time such that $g_T(v_2) \in W^s_g(v_3)$ or $W^u_g(v_3)$.  In order to 
make a true ideal rectangle we require that the leaves connecting these pairs alternate 
between stable and unstable.  Note that in figure 2, the leaf containing the first pair, 
$\{v, v_1\}$ is a stable leaf, but we can similarly start with an unstable leaf.  We now 
make the following definition:

\begin{defn} Let $v$, $v_1$, $v_2$, $v_3$ be vectors in $T^1\tilde M$ describing an ideal 
rectangle as indicated above.  Let $\tilde{F}_T$ be the restriction of the time $T$ frame 
flow map to the frames based at $v_2$.  Note that the choice of $v_1$ and $v_3$ uniquely 
determines this rectangle and define
\[h(v_1, v_3) = p(v_3, v) \circ p(g_T(v_2), v_3) \circ \tilde{F}_T \circ p(v_1, v_2) 
\circ p(v, v_1). \]
Let $\hat{H}_v$ be the closure of the group generated by all such $h(v_1, v_3)$. 
\end{defn}

\begin{figure}
\centering
\includegraphics[height=.45\textheight, width=.6\textwidth]{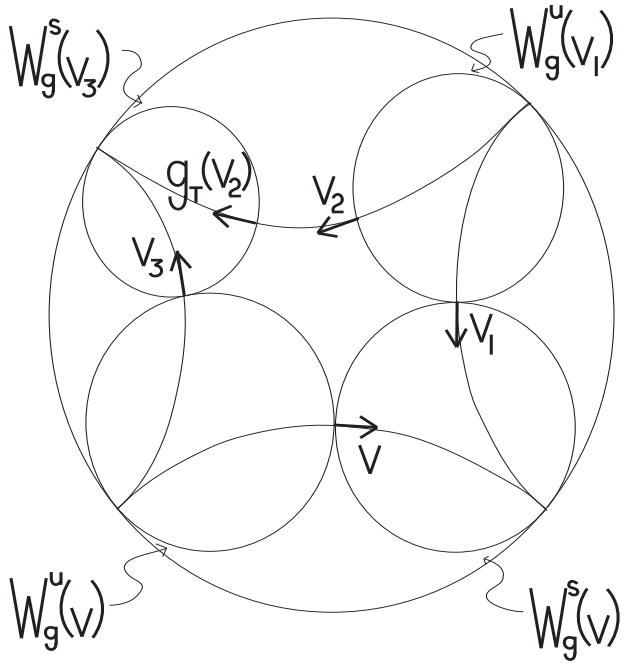}
\caption{general ideal rectangle}
\end{figure}

Note that $\hat{H}_v$ allows parallel translations along \emph{all} ideal rectangles 
based at $v$.  Furthermore, it is easy to see that a parallel translation around any 
`equilateral' ideal polygon as in the definition of $H_v$ can be broken up into a series 
of translations around the general ideal rectangles allowed in the definition of 
$\hat{H}_v$.  Thus, $\hat{H}_v \supseteq H_v$.  However, $\hat{H}_v$ preserves ergodic 
components, as frame flow certainly preserves ergodic components.  As the group $H_v$ is 
completely determined by the ergodic component of $n$-frames and $\hat{H}_v$ produces 
this same ergodic component, $\hat{H}_v \subseteq H_v$.  Therefore, $\hat{H}_v = H_v$.

Despite the equality of these groups, for our purposes there is a benefit to allowing 
this seemingly looser definition.  Each $(v_1, v_3) \in W^s_g(v) \times W^u_g(v)$ defines 
a rectangle used in $\hat{H}_v$.  This is in opposition to the case for $H_v$, where only 
a set of measure zero define allowed rectangles.  The advantage of this is the following. 
Since $M$ is negatively curved, the geodesic flow is ergodic and the unstable foliation 
$W^u_g$ is absolutely continuous.  Thus, there is a full measure set of $v_1 \in T^1M$ 
with dense negative time orbit and it must intersect the leaf $W^s_g(v)$ in a set of full 
conditional measure for almost all $v \in T^1M$ (see Appendix to \cite{Brin-ergodicity} 
Lemma 5.4).  Picking $v_1$ from this set will ensure the needed negative time recurrence 
of $\gamma_{v_2}$ since $\gamma_{-v_2} \to \gamma_{-v_1}$.  Likewise, for almost every $v 
\in T^1M$ a full conditional measure set of $v_3 \in W^u_g(v)$ will have the needed 
positive time recurrence under the geodesic flow.  We conclude that we can find a full 
measure set of $v$ having dense forward and backward orbits \emph{and} (using Fubini's 
theorem) with a full measure set of $W^s_g(v) \times W^u_g(v)$ possessing the desired  
dynamical properties for $(-v_1, v_3)$.  Therefore the desired recurrence properties are 
generic in the set of rectangles used to generate $\hat{H}_v$ for almost all $v$. The 
final fact needed to prove that the transitivity group preserves the distinguished vector 
fields is the following:

\begin{lemma}
The map $(v_1, v_3) \mapsto h(v_1, v_3)$ is continuous.
\end{lemma}

\begin{proof}
This Lemma follows from the fact that the frame flow admits a continuous foliation, which 
was proved by Brin \cite{Brin-toptrans}.

First, note that $p(v, v_1)$ and $p(v_3, v)$ are defined by leaves of the foliations for 
the frame flow.  The continuous dependence of these maps on $(v_1, v_3)$ follows 
precisely from the continuity of the leaves.

Second, as $(v_1, v_3)$ varies, the leaves $W^u_g(v_1)$ and $W^s_g(v_3)$ and the geodesic 
connecting $\gamma_{v_1}(-\infty)$ to $\gamma_{v_3}(\infty)$ all vary continuously.  Thus 
$v_2$, $T$ and $g_T(v_2)$ will vary continuously. Along with the argument of the previous 
paragraph, all this implies that the maps $p(v_1, v_2)$ and $p(g_T(v_2), v_3)$ depend 
continuously on $(v_1, v_3)$.  Also, $\tilde{F}_T$ will depend continuously on $(v_1, 
v_3)$ as the frame flow is continuous.

Since $h(v_1, v_3)$ is the composition of these maps, the Lemma is proved.
\end{proof}

Now we can prove the following result:

\begin{prop} \label{prop:preserves}
For almost all $v \in T^1M$, if $w_v(t)$ is a parallel field along $\gamma_v(t)$ making 
constant curvature $-a^2$ with the geodesic direction, then for every element $h$ in the 
transitivity group we have $K(h(w_v(0)), v) = -a^2$. \end{prop}

\begin{proof}
As discussed above, for almost all $v \in T^1M$ a full measure set of the $(v_1, v_3) \in 
W^s_g(v) \times W^u_g(v)$ give rectangles with the recurrence properties necessary for 
$h(v_1, v_3) \in \hat{H}(v)$ to map $w_v(0)$ to another distinguished vector field along 
$\gamma_v$.  In particular, this set of `nice' $(v_1, v_3)$ is dense in $W^s_g(v) \times 
W^u_g(v)$.  Since $h(v_1, v_3)$ depends continuously on $(v_1, v_3)$ and for a dense set 
of $(v_1, v_3)$, it preserves the distinguished fields, we have that all $h(v_1, v_3)$ 
preserve the distinguished fields.  Since $\hat{H}_v = H_v$ is generated by these 
elements, the transitivity group preserves the curvature $-a^2$ as desired.
\end{proof}

This result gives us the desired relationship between the transitivity group and the 
distinguished vector fields.  We can now apply the results of Brin-Gromov and 
Brin-Karcher and prove Theorem 1 easily.

\begin{mainthm}
Let $M$ be a compact, negatively curved manifold.  Suppose that along every geodesic in 
$M$ there exists a parallel vector field making sectional curvature $-a^2$ with the 
geodesic direction.  If $M$ is odd dimensional, or if $M$ is even dimensional and 
satisfies the sectional curvature pinching condition $-\Lambda^2 < K < -\lambda^2$ with 
$\lambda/\Lambda > .93$ then $M$ has constant negative curvature equal to $-a^2$.
\end{mainthm}

\begin{proof}
We showed in Proposition \ref{prop:preserves} that for almost all $v \in T^1M$ the 
sectional curvature $K(h(w_v(0)), v) = -a^2$ for all $h$ in the transitivity group.  In 
the setting of the theorem, the results of Brin-Gromov and Brin-Karcher tell us that the 
2-frame flow is ergodic.  In particluar, since the transitivity group gives the ergodic 
component for this flow, the transitivity group must act transitively on $v^\perp \subset 
T^1M$.  Thus, $K(\cdot, v)$ is identically $-a^2$.  Since this holds for almost all $v$ 
it holds for all $v$ by continuity of $K(\cdot, \cdot)$, and the theorem is proved.
\end{proof}


\section{Parallel fields and Jacobi fields} \label{sec:parallel}

In \cite{sphericalrank} a distinction is made between `weak' and `strong' rank.  The 
existence of \emph{parallel} fields making extremal curvature is called strong rank; the 
existence only of \emph{Jacobi} fields making extremal curvature is called weak rank.  A 
parallel field can be scaled (by a solution to the real variable version of the Jacobi 
equation where the standard derivative replaces the covariant derivative) to produce a 
Jacobi field.  Thus, a proof under the less stringent hypothesis of weak rank implies a 
proof for strong rank.  Hamenst\"{a}dt's is the sole result prior to this paper for weak 
rank.  She states her main theorem for parallel fields only, but she shows in Lemma 2.1 
that in negative curvature a Jacobi field making maximal curvature is a parallel field 
scaled by a function \cite{Hamenstadt}.  Essentially, she shows that Jacobi fields making 
maximal curvature grow at precisely the rate one finds for the constant curvature case.  
Connell accomplishes the same in \cite{Connell} Lemma 2.3.  This, together with some of 
the arguments below, shows that these Jacobi fields are in fact parallel fields scaled by 
an appropriate function.  Therefore, Corollary 2 is a weak rank result, needing only the 
Jacobi field hypothesis.

In this section we show that Jacobi fields making \emph{minimal} curvature with the 
geodesic direction are also scaled parallel fields.  This will justify the phrasing of 
Corollary 1 as a weak rank result.

First, note that we need only consider non-vanishing Jacobi fields; hence it will be 
enough to prove that stable and unstable Jacobi fields are scaled parallel fields.  
Stable Jacobi fields are those which have norm approaching zero as $t \to \infty$; 
unstable Jacobi fields have the same property in the negative time direction.  Suppose 
$J(t)$ is a stable Jacobi field along the geodesic $\gamma(t)$ making curvature $-a^2$ 
with the geodesic (take $a>0$ now), where $-a^2$ is the curvature minimum for the 
manifold (the modifications of what follows for unstable Jacobi fields are 
straightforward).  The Rauch Comparison Theorem (see \cite{doC} Chapt 10, Theorem 2.3) 
can be used to show that

\begin{equation} \label{eq:geq}
|J(t)| \geq |J(0)| e^{-at}.
\end{equation}

We would like to show that equality is achieved in (\ref{eq:geq}).  Write $J(t)=j(t)U(t)$ 
where $j(t)=|J(t)|$ and $U(t)$ is a unit vector field.  Then the Jacobi equation for $J$ 
reads:

\begin{equation} \label{eq:jacobi}
j''U + 2j'U' + jU'' + jR(\dot{\gamma}, U)\dot{\gamma} = 0
\end{equation}
where $j'$ denotes the standard derivative and $U'$ denotes covariant derivative.  Taking 
the inner product of (\ref{eq:jacobi}) with $U$ and noting that $\langle U'', U\rangle = 
- \langle U', U' \rangle$ we obtain

\begin{equation} \label{eq:simpjacobi}
j'' - j(\langle U', U'\rangle +a^2) = 0.
\end{equation}

We now know the following about the magnitude of $J$: $j \geq 0$ by definition, $\lim_{t 
\to \infty} j(t) = 0$ since $J$ is a stable Jacobi field, and $j'' \geq a^2 j$ by 
(\ref{eq:simpjacobi}).  These allow the following conclusion; its proof was shown to the 
author by Jeffrey Rauch:

\begin{lemma} \label{Rauchlemma}
Let $j$ be a non-negative, real valued function satifsying $j'' \geq a^2 j$ and $\lim_{t 
\to \infty} j(t) = 0$.  Then $j(t) \leq j(0) e^{-at}$.
\end{lemma}

\begin{proof}
We have that $a^2j-j'' \leq 0.$  On the interval $0 \leq t \leq R$ for $R\gg1$ define 
$g_R$ by $g_R(0)=j(0)$, $g_R(R)=j(R)$ and $a^2g_R-g_R'' = 0$.  Note that as $R \to 
\infty$, $g_R \to j(0)e^{-at}$.  We claim that $j \leq g_R$; the Lemma follows in the 
limit.

This claim is essentially the maximum principle.  First, $j \leq g_R$ holds at $0$ and 
$R$.  Now suppose $j-g_R$ has a positive maximum at $c \in (0, R)$.  Then $(j''-g_R'')(c) 
\leq 0$.  However, we know $a^2(j-g_R) - (j''-g_R'') \leq 0$, so a positive value of 
$j-g_R$ at $c$ together with a negative value of $j''-g_R''$ yields a contradiction.  
Therefore $j \leq g_R$ holds on all of $[0, R]$ as desired.
\end{proof}

This Lemma, together with equation (\ref{eq:geq}), tells us that $|J(t)| = 
|J(0)|e^{-at}$.  Examining equation (\ref{eq:simpjacobi}) we see that having the growth 
rate $e^{-at}$, as in the constant curvature $-a^2$ case, implies that $U' =0$, that is, 
$J$ is a scaled parallel field, as desired.


\section{The dynamical perspective} \label{sec:Lyap}

In this section we discuss how the results of Connell in \cite{Connell} can be adapted to 
prove Theorem 2 as a simple consequence of Corollary 1.  The necessary changes are for 
the most part cosmetic; the discussion here is included for completeness, but the author 
does not claim to have added anything of substance to Connell's work.  The notation below 
that has not already been assigned follows Connell's for ease of reference.

Recall that Lyapunov exponents are a tool for measuring long-term asympotic growth rates 
in dynamical systems (see Katok and Mendoza's Supplement to \cite{H-K} section S.2 for an 
exposition).  In the setting of the geodesic flow they can be defined as follows.  Let $v 
\in T^1M$ and $u \in v^{\perp}$.  Let $J_u(t)$ be the unstable Jacobi field along 
$\gamma_v$ with initial condition $J_u(0) = u$.  Then, the \emph{positive Lyapunov 
exponent at $v$ in the $u$-direction} is
\[\lambda_v^+(u) = \limsup_{t \to \infty} \frac{1}{t} log|J_u(t)|. \]
Define \[\lambda_v^+ = \max_{u \in v^{\perp}} \lambda_v^+(u). \]
This is the maximal Lyapunov exponent at $v$; the curvature bound $-a^2 \leq K$ (again, 
take $a>0$) implies that $\lambda_v^+ \leq a$.  Let
\[ \Omega = \{v \in T^1M : \lambda_v^+ = a\}. \]
We can now rephrase Theorem 2 more succinctly.

\begin{thm2}
Let $M$ be a compact manifold with sectional curvature $-a^2 \leq K < 0$.  Suppose that 
$\Omega$ has full measure with respect to a geodesic flow-invariant measure $\mu$ with 
full support.  If $M$ is odd dimensional, or if $M$ is even dimensional and satisfies the 
sectional curvature pinching condition $-a^2 \leq K < -\lambda^2$ with $\lambda/a > .93$ 
then $M$ is of constant curvature $-a^2$.
\end{thm2}

Connell shows in the upper rank case that the dynamical assumption implies the geometric 
one, that is, that the manifold in fact has higher rank, allowing the application of an 
appropriate rank rigidity theorem.  He first shows (\cite{Connell} Proposition 2.4) that 
along a closed geodesic $\lambda_v^+ = a$ implies the existence of an unstable Jacobi 
field making curvature $-a^2$ with the geodesic direction.  Essentially, if the Jacobi 
field giving rise to the Lyapunov exponent does not have this curvature, it will 
continually see non-extremal curvature a positive fraction of the time as it moves around 
the closed geodesic.  This contradicts the supposed value of the Lyapunov exponent.  The 
lower curvature bound version of the argument is exactly the same as that presented by 
Connell, with the proper inequalities reversed; also note that the work in section 
\ref{sec:parallel} of this paper gives the results analogous to Connell's Lemma 2.3 
necessary for the argument.

It is clear that if a dense set of geodesics have the distinguished Jacobi fields, then 
all geodesics will.  Since the velocity vectors for closed geodesics are dense in $T^1M$, 
Connell finishes his proof in section 3 of \cite{Connell} by showing that these vectors 
are all in $\Omega$ and using the argument of the previous paragraph.  Adapted to the 
setting of Theorem 2 the argument runs as follows.  If $w \in T^1M$ is tangent to a 
closed geodesic and $\lambda_w^+ < a$ the previous paragraph implies that any unstable 
Jacobi field along $\gamma_w$ must make curvature strictly greater than $-a^2$ for a 
positive amount of time.  By continuity, this will also be true of any unstable Jacobi 
field along a geodesic $\gamma_v$ in a sufficiently small neighborhood of $\gamma_w$ (in 
the Sasaki metric on $T^1M$).  The ergodic theorem implies that for a full measure set of 
$v \in T^1M$, $\gamma_v$ will spend a positive fraction of its life in this small 
neighborhood of the periodic geodesic $\gamma_w$; the positivity follows from the fact 
that $\mu$ has full support.  The intersection of this full measure set with the full 
measure set $\Omega$ thus contains vectors $v$ which have $\lambda_v^+ = a$ but spend a 
positive fraction of their life so close to $\gamma_w$ that no Jacobi fields along them 
can make the minimal curvature $-a^2$ with the geodesic direction during this fraction of 
the time.  In fact, since $\gamma_w$ is compact, so is the closure of this small 
neighborhood and therefore the curvature between these Jacobi fields and the geodesics, 
when in this neighborhood, can be bounded away from $-a^2$, i.e. $K(J_u, 
\dot{\gamma}_v)>c>-a^2$ for a fixed $c$.  Having this curvature bound a positive fraction 
of the time contradicts $\lambda_v^+ = a$; therefore all closed geodesics must lie in 
$\Omega$ and the argument is complete.

Again, this version of the argument, relevant for the lower curvature bound situation, is 
the same as that presented by Connell with the proper inequalities reversed.  Thus, the 
dynamical assumption implies the geometric assumption of Corollary 1 and Theorem 2 
follows.  Note that for these arguments the extremality of the distinguished curvature is 
essential; a result that parallels Theorem 1 in allowing non-extremal distinguished 
curvature cannot be hoped for.


\section{Conclusion} \label{sec:conclusion}

We conclude with a few remarks on our results in the context of the other rank rigidity 
theorems.  As noted above, Corollary 2 treats the case dealt with by Hamenst\"{a}dt's 
hyperbolic rank rigidity theorem, strong upper hyperbolic rank.  Unlike Hamenst\"{a}dt's 
result, the result presented in this paper is limited by the curvature pinching condition 
in even dimension.  However, this paper's proof is shorter, and has the advantage of 
telling us which symmetric space $\tilde{M}$ is.  Corollary 1 is strong lower hyperbolic 
rank rigidity, and this result is the first positive result for lower rank rigidity of 
any sort.  Counterexamples to lower rank rigidity in other curvature settings are known; 
\cite{sphericalrank} presents an overview, together with counterexamples to weak upper 
and lower spherical rank rigidity.  Counterexamples to lower Euclidean rank rigidity are 
given in \cite{Spatzier-Strake}. When the value $-a^2$ is not extremal we have a result 
of a different type than previous rank rigidity results.  Our results also show that 
spaces of constant negative curvature can not be deformed while maintaining the 
distinguished parallel vector fields along all geodesics, or while maintaining extremal 
Lyapunov exponents at a full measure set of $T^1M$, except by scaling the metric.

Note that in even dimension a result as extensive as our odd dimensional result cannot be 
hoped for.  Since parallel translation preserves the complex structure on a K\"{a}hler 
manifold the 2-frame flow will not be ergodic (see \cite{Brin-Gromov} for some results on 
unitary frame bundles).  These known counterexamples to ergodic frame flow are excluded 
by requiring $-1<K<-1/4$, leading Brin to conjecture that strict 1/4-pinching implies 
that the frame flow is ergodic (\cite{Brin-survey} Conjecture 2.6).  A positive answer to 
this conjecture, or any extended results for ergodicity of the 2-frame flow in negative 
curvature would extend the results on rank rigidity presented here correspondingly, using 
the same proof as presented above.  One still hopes that lower hyperbolic rank rigidity 
(in the sense that higher rank implies the space is locally symmetric) could be true 
without any curvature pinching in even dimensions, perhaps even without the restriction 
$K < 0$, but such a result would call for a completely different method of proof from 
that presented here.

\end{document}